\documentclass[smallextended]{article}       
\usepackage[utf8]{inputenc} 
\usepackage{amssymb,amsmath} 
\usepackage{graphicx} 
\usepackage{array} 
\usepackage{eurosym} 
\usepackage{xcolor} 
\usepackage{booktabs}
\usepackage{epstopdf}
\usepackage{pifont} 
\usepackage{algorithm} 
\usepackage{algpseudocode}
\usepackage{multicol}
\usepackage{multirow}
\usepackage{subcaption}
\usepackage{changepage}
\usepackage{hyperref}
\usepackage{adjustbox}
\usepackage{graphicx}
\usepackage{multirow}
\usepackage{amsmath} 
\usepackage[affil-it]{authblk}

\begin{document}

\title{A novel bidding method for combined heat and power units in district heating systems}


\author[1]{Ignacio~Blanco}
\author[2]{Anders N. Andersen}
\author[1,*]{Daniela Guericke}
\author[1]{Henrik~Madsen}

\affil[1]{\small Technical University of Denmark, Department for Applied Mathematics and Computer Science, Richard Petersens Plads, 2800 Kgs. Lyngby, Denmark}
\affil[2]{\small EMD International A/S, Niels Jernesvej 10, 9220 Aalborg \O, Denmark}
\affil[*]{\small Corresponding author: Daniela Guericke, dngk@dtu.dk}

\maketitle

\begin{abstract}
We propose a bidding method for the participation of combined heat and power (CHP) units in the day-ahead electricity market. More specifically, we consider a district heating system where heat can be produced by CHP units or heat-only units, e.g., gas or wood chip boilers.  We use a mixed-integer linear program to determine the optimal operation of the portfolio of production units and storages on a daily basis. Based on the optimal production of subsets of units, we can derive the bidding prices and amounts of electricity offered by the CHP units for the day-ahead market. The novelty about our approach is that the prices are derived by iteratively replacing the production of heat-only units through CHP production.  This results in an algorithm with a robust bidding strategy that does not increase the system costs even if the bids are not won. We analyze our method on a small realistic test case to illustrate our method and compare it with other bidding strategies from literature, which consider CHP units individually. The analysis shows that considering a portfolio of units in a district heating system and determining bids based on replacement of heat production of other units leads to better results.
\end{abstract}

\section{Introduction}
The global target of reducing CO2-emissions from fossil fuels has required several countries, especially in the European Union, to consider efficient district heating and cooling systems as a key role in its CO2-emissions reduction strategy \cite{studyone27:online}. Since it is assumed that fossil fuels will be mostly replaced by intermittent renewable energy sources, a higher share of district heating and cooling systems can facilitate the integration of these intermittent energy sources in the energy mix \cite{CONNOLLY2014475}, contributing
to balance the grid by the use of heat pumps, electric boilers, thermal storage or flexible CHP production. The efficiency of these systems has been demonstrated already in countries like Denmark and Sweden. In Denmark around 64\% of the households are connected to district heating networks for space heating and domestic hot water \cite{euroheat}. Nowadays, the total heat consumption in all district heating networks in Denmark is close to 130 petajoule (PJ) from which more than 65\% are produced by combined heat and power (CHP) plants \cite{regulati50:online}. The integration of renewable and intermittent energy sources (e.g. wind and photovoltaic) in liberalized electricity markets such as Nordpool \cite{NordPool52:online} lead to historical low and highly volatile electricity prices \cite{morales2013integrating,Jaehnert2014}. This results in a larger difficulty for CHP units to be scheduled and thus obtain profits from the electricity market \cite{Jaehnert2014}. Consequently, they are being replaced by other heat production units such as heat pumps that take advantage of low electricity prices in periods where the mix of renewable energy production is high \cite{ostergaard2016booster,nyborg2015heat}. This means many district heating companies are no longer operating the system only with CHP units but in combination with other heat production units such as gas boilers, heat pumps, electric boilers or wood chip boilers. In this work, we propose a method that optimizes the power production bids of CHP units in district heating systems in the day-ahead market, making these units more competitive and scheduled in more hours. In addition, the approach provides flexibility to the power system by activating or deactivating the CHP units when required by the transmission system operator (TSO), which is reflected in the market prices. 

The optimal operation of standalone CHP units has been extensively considered in literature \cite{dimoulkas2015probabilistic,illerhaus1999optimal,schaumburg2006partial}. In these publications, the authors use mixed-integer linear programming to define the technical constraints of the CHP unit in combination with a thermal storage to maximize its profits in a liberalized electricity market framework.

The use of mathematical optimization to operate a portfolio of heat production units has been proposed in several publications. The most relevant publications with respect to our approach are reviewed in the following. First, the work of \cite{fang2016optimization} shows that adding a receding horizon to optimize the operation of a thermal storage yields better results in terms of costs. In \cite{li2016optimal} the authors consider an integrated power and heat system. The goal of their approach is to coordinate the heat and power production to accommodate the intermittent generation of a wind farm in the district heating network by using a heat pump when the wind power production is high. Following the same principle of producing heat  when the electricity price is low and produce power and heat with the CHP when the electricity price is high, the authors in \cite{nielsen2016economic} evaluate the scheduling of different heat production units under electricity price and heat demand uncertainty with two-stage stochastic programming. Finally, the approach proposed in \cite{ito2017electricity} integrates many different heat and cooling production units and use optimization to operate them as a portfolio. The system offers capacity to the electricity market to accommodate power fluctuations by the interplay of the different units. Apart from scientific literature, there exist commercial software tools, such as \textit{energyPRO} \cite{EMDInter28:online} that schedule the production in integrated district heating systems. 

The above mentioned literature considers the optimal production in district heating systems including CHP units, but neglects to determine bids that the operator should present to the electricity market. CHP units are usually more expensive to operate than other heat production units. Therefore, it has to be ensured that no economic losses due to poor bidding strategies will occur while operating the system. A typical liberalized electricity market works on a short term basis. An auction takes place on the day before the energy is delivered, the so-called day-ahead market \cite{Scharff2014}. Further markets exist such as the reserve capacity market, intra-day market and balancing market \cite{Scharff2014,experience2005lessons,sioshansi2011competitive}. The price of electricity from one day to another is partly unknown and more unpredictable as more intermittent renewable energy technologies get into the system. Despite that methods to predict electricity prices as a function of intermittent energy sources have been proposed (e.g. \cite{jonsson2013forecasting}), the optimal operation of  CHP units cannot be determined one day in advance when the power production bids have to be submitted, since optimal scheduling is depending on electricity price forecasts.

To approach this problem, we propose a method that considers the optimal operation as well as the optimal bidding for CHP units in district heating systems. Similarly, several bidding methods for both thermal power generation and CHP units have been proposed in literature. In the following, we review those methods, which are also used for comparison in the numerical results.

First, we review bidding methods used for thermal power generation units. Although these consider the power generation unit as a standalone production unit, these methods could be used by district heating operators to determine bids for the CHP plant without taking the other units into account.  One strategy is proposed by \cite{conejo2002price}. Their bidding procedure consists of creating bounds on the uncertainty given by forecasted prices and use these bounds to generate offers. They ensure profitability in their bids by offering power volumes at two different prices, i.e., the defined lower bound and upper bound of their price forecast. Since the vast majority of possible realizations of the electricity price lie within the defined bounds, they protect their offers against uncertainty. The second method analyzed is the method presented in \cite{rodriguez2004bidding}. This method generates different confidence intervals for electricity prices. By solving the model for each of these intervals, bidding curves can be created.

Second, we mention bidding strategies where the authors consider a single CHP unit with a thermal storage. The method proposed in \cite{schulz2016optimal} determines the optimal production of the CHP unit. The bidding price is the price forecast, which is the same price used to determined the power production. In \cite{dimoulkas2014constructing} the authors construct different bidding curves for a CHP unit based on price scenarios and calculate for each offering period the power production vs. the electricity prices. These curves are submitted to the day-ahead market.  The mentioned publications do not take  advantage of the different heat production technologies that are connected to the district heating network. The potential of the system can be increased if these are optimized jointly as a portfolio.

Finally, in \cite{ravn2004modelling} the authors propose a bidding strategy for CHP units that takes into account other heat units to define the heat production costs. Therefore, the bidding prices generated by this method are the same as those generated in ours. However, the way the volume and the bidding hours are calculated is different, because \cite{ravn2004modelling} use a piece-wise linear function to activate different volumes of power at different prices according to the price forecast.

In this publication we introduce a novel bidding method for the participation of combined heat and power (CHP) units that operated in a portfolio with other heat production units in the day-ahead electricity market based on mixed-integer linear programming. The main contributions of our work are: 
\begin{itemize}
	\item Our method operates the heat production units as a portfolio and takes the costs of all units in the system into account to define the optimal bidding prices. We show that this is more beneficial for district heating providers than using state-of-the-art bidding methods for standalone CHP or thermal units. 
	\item The bids generated by our method consider the cost of producing heat to define bids, which protects against the uncertainty of electricity prices, i.e., lost bids do not increase the operational cost. 
	\item We develop an iterative process to generate our bids based on replacement of heat production by heat-only units that results in a higher number of offers compared to other bidding methods.  
	\item We reimplemented the above mentioned methods for operational optimization and bidding \cite{conejo2002price,rodriguez2004bidding,schulz2016optimal,dimoulkas2014constructing,ravn2004modelling} to compare them to our bidding strategy in a case study.
	\item To improve the daily operation of a heat storage system, we analyze different lengths of receding horizon.
	
\end{itemize}

The remainder of the paper is organized as follows. In Section \ref{sec:model} the operational planning problem is explained in detail. In Section \ref{sec:biddingmethod}, we develop step by step the bidding method proposed in this work. Section \ref{casestudies} presents the case study used to run our experiments. In Section \ref{sec:illustrativeexample}, an illustrative example shows how the bidding method works. Further simulations and numerical results are presented in Section \ref{sec:results}. Finally, we summarize our work in Section \ref{conclusion}. 

\section{Operational planning model}\label{sec:model}
We start by introducing the mixed-integer linear program (MILP) to schedule the optimal operation of a portfolio of heat production units in a district heating system. For an overview of the nomenclature, we refer to Table \ref{tab:nomenclature}. 

\begin{table}
	\footnotesize
	\centering
	\caption{Nomenclature}
	\begin{tabular}{p{0.16\textwidth}p{0.7\textwidth}}\toprule
		\multicolumn{2}{l}{\textbf{Sets}}\\\midrule
		$\mathcal{T}$ & Set of time periods $t$\\
		$\mathcal{U}$ & Set of heat production units $u$\\
		$\mathcal{U}^{\text{CHP}} \subset\mathcal{U} $& Subset of combined heat and power production units\\
		$\mathcal{U}^{\text{H}} \subset\mathcal{U}$ & Subset of heat-only production units \\\midrule
		\multicolumn{2}{l}{\textbf{Parameter}}  \\\midrule
		$C^{\text{}}_{u}$ & Cost for producing heat with unit $u \in \mathcal{U}$  [DKK/MWh-heat]\\
		$\overline{Q}_{u}$ &Max. heat production for unit $u \in \mathcal{U}$ [MWh-heat]\\
		$\underline{Q}_{u,t}$ &Min. heat production for unit $u \in \mathcal{U}, t \in \mathcal{T}$   [MWh-heat]\\
		$A^{\text{DH}}_{u}$ & Binary parameter: 1, if unit $u \in \mathcal{U}$ is connected to the district heating system, 0, otherwise\\
		$A^{\text{S}}_{u}$ & Binary parameter: 1, if unit $u \in \mathcal{U}$ is connected to the thermal storage, 0, otherwise\\
		$\varphi_{u}$ &Heat-to-power ratio for unit $u \in \mathcal{U}^{\text{CHP}}$ [$\text{MWh-heat}/\text{MWh-el}$]\\
		$\overline{P}_{u}$& Max. power production for unit $u \in \mathcal{U}^{\text{CHP}}$ [MWh-el] \\
		$\underline{P}_{u}$ & Min. power production for unit $u \in \mathcal{U}^{\text{CHP}}$ [MWh-el]\\
		$S^{\text{F}}$& Maximum heat flow from the storage to the district heating network [MWh-heat] per period\\
		$s_{0}$ & Initial storage level [MWh-heat]\\
		$\overline{S^{\text{}}}$& Maximum heat storage level [MWh-heat]\\
		$\underline{S^{\text{}}}$& Minimum heat storage level [MWh-heat]\\
		${\lambda}_{t}$ &Electricity price forecast for time period $t \in \mathcal{T}$ [DKK/MWh-el]\\
		$D_{t}$ &Heat demand for time period $t \in \mathcal{T}$ [MWh-heat]\\\midrule
		\multicolumn{2}{l}{\textbf{Variables}}  \\\midrule
		$q_{u,t} \in \mathbb{R}^+_{0}$ &Heat production of heat unit $u \in \mathcal{U}$ in period $t \in \mathcal{T}$ [MWh-heat]  \\  
		$q_{u,t}^{\text{DH}} \in \mathbb {R}^+_{0}$ & Heat production of unit $u \in \mathcal{U}$ inserted to the grid in period $t \in \mathcal{T}$ [MWh-heat] \\  
		$q_{u,t}^{\text{S}} \in \mathbb{R}^+_{0}$ & Heat production of unit $u \in \mathcal{U}$ inserted to the storage in period $t \in \mathcal{T}$ [MWh-heat]  \\
		$p_{u,t} \in \mathbb{R}^+_{0}$ & Power production of unit $u \in \mathcal{U}^{\text{CHP}}$ in period $t \in \mathcal{T}$ [MWh-el] \\
		$s_{t} \in \mathbb{R}^+_{0}$ & Thermal storage level at time period $t \in \mathcal{T}$ [MWh-heat]  \\
		$s^{\text{OUT}}_{t} \in \mathbb{R}^+_{0}$ & Heat flowing from the storage to the district heating in period $t \in \mathcal{T}$ [MWh-heat]  \\
		$x_{u,t} \in \{0,1\}$ & Binary variable: 1, if unit $u \in \mathcal{U}^{\text{CHP}}$ is on in period $t \in \mathcal{T}$ and 0, otherwise\\\bottomrule
	\end{tabular}
	\label{tab:nomenclature}
\end{table}

The set of considered production units is denoted by $\mathcal{U}$, which is comprised of CHP units $\mathcal{U}^{\text{CHP}}$ (producing heat and power simultaneously) and heat-only units $\mathcal{U}^{\text{H}}$. We consider the production over a time horizon of $\mathcal{T}$ periods, where each period is one hour due to the hourly bidding periods in the day-ahead market.  Each unit $u \in \mathcal{U}$ has a maximum heat production per hour $\overline{Q}_{u}$ and production cost of $C^{\text{}}_{u}$ per MWh heat. The CHP units $u \in \mathcal{U}^{\text{CHP}}$ have further a restriction on the minimum and maximum electricity production per hour denoted by $\underline{P}_{u}$ and $\overline{P}_{u}$, respectively. The heat-to-power ratio for the CHP units is given by $\varphi_{u}$. 
The thermal storage of the system has a minimum ($\underline{S^{\text{}}}$) and maximum level ($\overline{S^{\text{}}}$) and the outflow per period is limited to $S^{\text{F}}$ (all in MWh heat).  The binary parameters $A^{\text{DH}}_{u}$  and $A^{\text{S}}_{u}$  determine whether a heat unit $u$ is connected to the district heating network and/or the thermal storage, respectively. 
The heat demand in the district  heating network for each periods is given by $D_{t}$ and the forecasted electricity price for each hour by ${\lambda}_{t}$.

The model determines the optimal heat production for all units $u \in \mathcal{U}$ in each period in variables $q_{u,t}$ where the production can either go directly to the district heating network ($q_{u,t}^{\text{DH}}$) or to the thermal storage ($q_{u,t}^{\text{S}}$). The power production of the CHP units $u \in \mathcal{U}^{\text{CHP}}$ is modelled by variables $p_{u,t}$ and their status (on/off)  by the binary variable $x_{u,t}$. The variables $s_{t}$ and $s^{\text{OUT}}_{t}$ represent the storage level and storage outflow in each period $t$.
Note that the number of periods $\mathcal{T}$ is not limited to 24 hours, because it can be profitable to already consider the production of future days to operate the thermal storage more efficiently. The bidding method proposed in the next section is used to create hourly bids for the day-ahead market on a daily basis, thus, we shift the planning period by 24 hours after each run in a receding horizon approach.

{
	\begin{subequations}
		\allowdisplaybreaks
		The objective function \eqref{objectivefunction1} minimizes the cost of producing heat by units $\mathcal{U}$, while taking the expected income from the electricity market into account.
		\begin{align}
		\underset{}{\text{min}} & \hspace{0.2cm}\sum_{t \in T} \sum_{u \in \mathcal{U}^{\text{CHP}}}\big(C_{u}\varphi_{u}-\lambda_{t} \big)p_{u,t}+\sum_{t \in T} \sum_{u \in \mathcal{U}^{\text{H}}} \big( C_{u}q_{u,t} \big)
		\label{objectivefunction1}\end{align}
		The production and flow of heat is modeled in constraints \eqref{Maxheatproduction} to \eqref{restricts}. The heat production capacity of all units is limited by constraint \eqref{Maxheatproduction}. The lower bound of the heat production is normally $\underline{Q}_{u,t}=0$ for all units $u \in \mathcal{U}$ and time periods $t \in \mathcal{T}$, but it is a necessary restriction that we will use in some parts of our algorithm which is described in more detail in Section \ref{sec:biddingmethod}. The produced heat is either used in the district heating system or flows to the thermal storage \eqref{storagedistrictheating}. Whether the unit is connected to the district heating, the storage or both is determined by constraints \eqref{restrictdh} and \eqref{restricts}.
		\begin{align}
		& \underline{Q}_{u,t} \leq q_{u,t} \leq \overline{Q_{u}} \label{Maxheatproduction} && \forall t \in \mathcal{T}, \forall u \in \mathcal{U}  \\
		&  q_{u,t} = q^{\text{DH}}_{u,t}+q^{\text{S}}_{u,t}  \label{storagedistrictheating} && \forall t \in \mathcal{T}, \forall u \in \mathcal{U}\\
		&  q^{\text{DH}}_{u,t} \leq \overline{Q_{u}}A^{\text{DH}}_u  \label{restrictdh} && \forall t \in \mathcal{T}, \forall u \in \mathcal{U}\\
		&  q^{\text{S}}_{u,t}\leq \overline{Q_{u}}A^{\text{S}}_u \label{restricts} && \forall t \in \mathcal{T}, \forall u \in \mathcal{U}
		\end{align}
		The power and heat production of a CHP unit is connected in constraints \eqref{Heattopower} with the corresponding heat-to-power ratio. Constraints \eqref{maxminpower} set the minimum and maximum power production of the CHP plant, if the CHP is on ($x_{u,t}=1$) and to 0 otherwise ($x_{u,t}=0$).
		\begin{align}
		& q_{u,t} = \varphi_{u}p_{u,t} \label{Heattopower}  &&  \forall t \in \mathcal{T}, \forall u \in \mathcal{U}^{\text{CHP}}   \\
		&  \underline{P}_{u}x_{u,t} \leq p_{u,t} \leq \overline{P}_{u}x_{u,t} \label{maxminpower}&&  \forall t \in \mathcal{T},  \forall u \in \mathcal{U}^{\text{CHP}}  
		\end{align}
		The heat storage is modeled through constraints \eqref{eq2:thermal_level}-\eqref{eq2:thermal_end}. Constraints \eqref{eq2:thermal_level} determine the storage level in each period while the capacity of storage is ensured by constraints \eqref{eq2:thermal_limits}. The inflow to and outflow from the storage is limited by constraints \eqref{eq2:thermal_maxoutflow} and \eqref{eq2:thermal_maxinflow}, respectively. We ensure that the initial storage level is reached again at the end of the receding planning horizon to avoid emptying the storage every day \eqref{eq2:thermal_end}. In Section \ref{sec:results}, we further investigate how many days of receding horizon should be considered.
		\begin{align}
		& s_{t} =s_{t-1}+\sum_{u \in \mathcal{U}_{}}q^{\text{S}}_{u,t}-s^{\text{OUT}}_{t}  && \forall t \in \mathcal{T}
		\label{eq2:thermal_level}\\
		&\underline{S}  \leq s_{t} \leq \overline{S}  && \forall t \in \mathcal{T} \label{eq2:thermal_limits}\\
		&\sum_{u \in \mathcal{U}_{}}q^{\text{S}}_{u,t}  \leq S^{\text{F}}  && \forall t \in \mathcal{T} \label{eq2:thermal_maxoutflow}\\
		&s^{\text{OUT}}_{t}  \leq S^{\text{F}}  && \forall t \in \mathcal{T} \label{eq2:thermal_maxinflow}\\
		&s_{|\mathcal{T}|}   \geq s_{0}  && \label{eq2:thermal_end} 
		\end{align}
		Finally, the heat balance of the system is ensured by equation \eqref{eq2:Heatbalance}, i.e., the heat output to the district heating network must match the heat demand in each period. We assume that the demand data is adjusted to take heat losses into account. 
		\begin{align}
		& D_{t} = \sum_{u \in \mathcal{U}}q^{\text{DH}}_{u,t}+s^{\text{OUT}}_{t}  && \forall t \in \mathcal{T} \label{eq2:Heatbalance}
		\end{align}
	\end{subequations}
The above stated MILP \eqref{objectivefunction1}-\eqref{eq2:Heatbalance} is the basis for our bidding method described in the next section.
}

\section{Bidding method}\label{sec:biddingmethod}
In this section we present our bidding method, named Heat Unit Replacement Bidding (HURB) method, for district heating operators. We assume that the district heating operator has enough capacity to cover the heat demand in all periods by just using heat-only units, i.e., units that are operated independently from won offers in the electricity markets. This is a reasonable assumption, because this is the case for almost all district heating systems in Europe. In particular, because not having enough capacity, such that they can operate independently from the electricity market, introduces a risk of not covering the heat demand, which is the primary goal of district heating providers.

The HURB method uses the MILP presented in Section \ref{sec:model} to optimize the heat production in each time period by using the units in $\mathcal{U}$. 
To obtain the hourly bidding prices and power amounts offered to the day-ahead from this solution, we proceed as follows. The idea behind our bidding method is to incentivize the CHP units to place as many bids as needed to replace the heat production  from heat-only units, i.e., selling power to the electricity market would lower the overall cost compared to just using heat-only units. This stems from the observation that CHP units are usually more expensive to operate than heat-only units. Thus, we would normally avoid producing with the CHP units if it is not necessary. However, if we have the chance to receive a sufficiently high price from the electricity market, such that the net costs (operational cost minus profits) of the CHP units drop below  production costs of some or all of the heat-only units, it is more profitable to produce heat with the CHP units.

\begin{table}
	\footnotesize
	\caption{Symbols in Algorithm \ref{alg:biddingmethod}}
	\begin{tabular}{lp{0.8\textwidth}}\toprule
		$i$ & Iteration counter\\
		$\mathcal{H}$ & Heat-only units in descending order of operational cost \\
		$\mathcal{U}^{\text{H}}_i$ & Subset of heat-only units to be considered in iteration $i$\\
		$\mathcal{O}$ & Set of offers\\
		$b_u$ & Bidding price for unit $u \in \mathcal{U}^{\text{CHP}}$\\
		$a_{u,t}$ & Bidding amount for unit $u \in \mathcal{U}^{\text{CHP}}$ in period $t \in \mathcal{T}$\\
		$k_{u,t}$ & Bidding amount for unit $u \in \mathcal{U}^{\text{CHP}}$ in period $t \in \mathcal{T}$ cummulated over the iterations\\
		$(u,t,b_u,a_{u,t})$ & Tuple representing an offer in set $\mathcal{O}$ valid for unit  $u \in \mathcal{U}^{\text{CHP}}$ in period $t \in \mathcal{T}$ with bidding price $b_u$ and amount $a_{u,t}$\\\bottomrule
	\end{tabular}
\label{tab:symbols}
\end{table}

\begin{algorithm}[t]
		\footnotesize
	\caption{Heat Unit Replacement Bidding (HURB)}
	\begin{algorithmic}
		\State Set $\lambda_t=0 \ \forall t \in T$  and $\underline{Q}_{u,t}=0 \ \forall u \in U, t \in T$
		\State Solve \eqref{objectivefunction1}-\eqref{eq2:Heatbalance} and store optimal heat production $q^*_{u,t}$
		\State Set $\lambda_t$ to the electricity price forecast
		\State $\mathcal{H}$ = sort($\mathcal{U}^H$) \Comment{Order heat-only units by descending $C_u$}
		\State$\mathcal{O} \gets \emptyset$ \Comment{Initialize set of offers  }
		\State $k_{u,t} \gets 0$ \Comment{Variable storing power production}
		\State $i\gets0$, $\mathcal{U}^{\text{H}}_i \gets \mathcal{U}^{\text{H}}$ \Comment{Initialize iteration counter and heat-only units}
		\For{each $h \in \mathcal{H}$}
		\State $i \gets i+1$
		\State $\mathcal{U}^{\text{H}}_i = \mathcal{U}^{\text{H}}_{i-1} \backslash \lbrace h \rbrace$ \Comment{Remove next heat-only unit}
		\State Set $\underline{Q}_{u,t} = q^*_{u,t} \forall u \in \mathcal{U}^{\text{H}}_i$
		\State Solve \eqref{objectivefunction1}-\eqref{eq2:Heatbalance} with $\mathcal{U} = \mathcal{U}^{\text{CHP}} \cup \mathcal{U}^{\text{H}}_i $ \Comment{Optimize }
		\State Get optimal solution $p^*_{u,t}$  
		\For{each $t \in \lbrace 1, \ldots 24 \rbrace$}
		\For {each $u \in \mathcal{U}^{\text{CHP}}$}
		\State $a_{u,t} \gets p^*_{u,t} - k_{u,t}$\Comment{Set amount}
		\If{ $a_{u,t} > 0$} \Comment{If new CHP production}
		\State $b^{}_{u} \gets  \big(C_{u}-C^{}_{h}\big)\cdot \varphi_{u}$ \Comment{Set bidding price}
		\State $\mathcal{O} \gets \mathcal{O} \cup \big(u,t,b_{u},a_{u,t}\big)$ \Comment{Add new offer}
		\EndIf
		\State $k_{u,t} \gets k_{u,t} + a_{u,t}$ \Comment{Cumulate power }
		\EndFor
		\EndFor
		\EndFor
		\State Return $\mathcal{O}$
	\end{algorithmic}
	\label{alg:biddingmethod}
\end{algorithm}

Based on this idea, we propose an iterative approach to determine when which heat-only unit will be replaced by CHP production. The outline  is shown in Algorithm \ref{alg:biddingmethod} and the symbols are defined in Table \ref{tab:symbols}. First, we obtain the optimal heat production for all units without consideration of participation in the electricity market. This means we solve the MILP \eqref{objectivefunction1} to \eqref{eq2:Heatbalance}  with $\lambda_t = 0 \ \forall t \in T$ and  no restriction on the minimum production, i.e., $\underline{Q}_{u,t} = 0 \ \forall u \in \mathcal{U}, t \in \mathcal{T}$ (see constraint \eqref{Maxheatproduction}). The optimal heat production of all heat-only units $u \in \mathcal{U}^{\text{H}}$ is stored in $q^*_{u,t}$ for later use. 
Then we sort the heat-only units in descending order of production costs per MW heat. From this we obtain the set $\mathcal{H} = \lbrace h_1, ..., h_n \rbrace$ where $\mathcal{H}$ is the ordered version of set $ \mathcal{U}^{\text{H}}$, i.e.,  $h_j \in \mathcal{U}^{\text{H}}$, $n = |\mathcal{U}^{\text{H}}|$ and $C_{h_j} > C_{h_k} \forall j,k \in \lbrace 1, \ldots, n \rbrace$ if $j<k$. Afterward, we iterate over the elements in set $\mathcal{H}$. In each iteration $i$ we remove the next most expensive unit $h$ from the set $\mathcal{U}^{\text{H}}_{i-1}$ resulting in a subset $\mathcal{U}^{\text{H}}_i = \mathcal{U}^{\text{H}}_{i-1} \backslash \lbrace h \rbrace$. For the remaining heat units $u \in \mathcal{U}^{\text{H}}_i$, we restrict the heat production with a lower bound $\underline{Q}_{u,t} = q^*_{u,t} $ to ensure that only  production of unit $h$ is replaced.  Then we determine the optimal production for this new subset of heat-only units and all CHP units $\mathcal{U}^{\text{CHP}}$ by solving the MILP \eqref{objectivefunction1} to \eqref{eq2:Heatbalance} for $\mathcal{U} =  \mathcal{U}^{\text{H}}_i\cup  \mathcal{U}^{\text{CHP}}$. From the solution we can determine the periods in which the heat-only unit removed in this iteration is replaced by the CHP units. This can be derived from the power production $p^*_{u,t}$ of the CHP units in the solution. The amount offered to the market $a_{u,t}$ is the additional power production added in this iteration. If the CHP unit is not producing, i.e., $p^*_{u,t}=0$, no offer is created.

The bidding price $b_{u}$ for offering the power produced by CHP unit $u$ is based on the production cost of the CHP unit $u$ and the replaced heat-only unit $h$ as shown in equation \eqref{biddingprice}.  
\begin{equation}
b^{}_{u} = \big(C_{u}-C^{}_{h}\big)\cdot \varphi_{u} \qquad \forall u \in \mathcal{U}^{\text{CHP}}
\label{biddingprice}
\end{equation}
The bidding price can be interpreted as follows. If we replace production by heat-only unit $h$ with CHP unit $u$, we have to pay the production cost of the CHP unit $u$ but we also save the production cost of heat-only unit $h$. To be profitable, we need to get at least the remaining amount $(C_{u}-C^{}_{h})$ from the market. This is our bidding price. At this price we are indifferent about whether we produce with the heat-only or with the CHP unit, because the costs for the production will be the same. By setting the bidding price in this way, we will not increase the cost even if we are not dispatched, because we have to produce the heat anyway. But if the market price is above our bidding price, we will make profits and lower our overall operational costs. 

We repeat this process by iteratively removing further heat-only units from the optimization problem until only the CHP units are left. During the iterations we collect the offers for the next day (i.e., periods $t \in \lbrace 1, \ldots 24 \rbrace$ although the total number of periods can be larger) in set $\mathcal{O}$. Each offer $o \in \mathcal{O}$ consist of a tuple $(u,t,b_{u}, a_{u,t})$ stating the CHP unit $u$ it corresponds to, the hour $t$ of the next day it is valid for, the bidding price $b_{u}$ and amount of power $a_{u,t}$. Offers created in the process are added to this set. No previous offers are removed, because they are still valid. After the algorithm is finished, the operator can determine the bids from set $\mathcal{O}$ and present them to the market.

\section{Case study}\label{casestudies} 
We use the following case study to demonstrate the HURB method in an illustrative example in Section \ref{sec:illustrativeexample} and compare the numerical results to other bidding methods in Section \ref{sec:results}.

\begin{figure}
	\centering
	\includegraphics[width=0.9\textwidth]{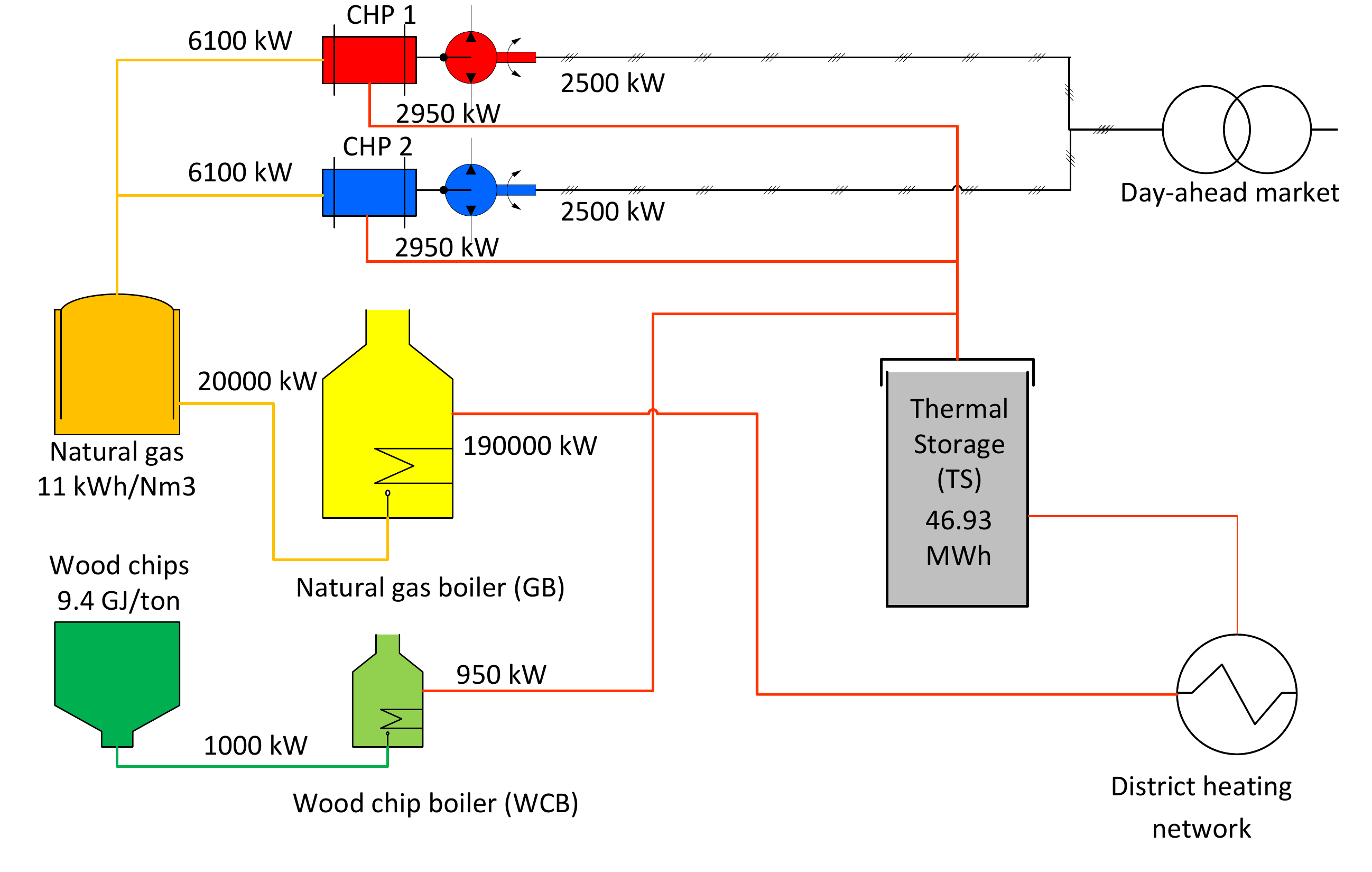}
	\caption{Flowchart of the system}
	\label{WCBGB2CHPS_Flowchart}
\end{figure}

We consider a small district heating system as shown in Figure \ref{WCBGB2CHPS_Flowchart}. The system includes two CHP production units (CHP1 and CHP2) in the form of gas engines and two heat-only units, more specifically a gas boiler (GB) and wood chip boiler (WCB). Furthermore, the system includes a thermal storage (TS) that can deliver heat to the district heating network. The CHP units and the WCB are connected to the thermal storage and not directly feeding to the district heating network. The GB is not connected to the thermal storage, but directly to the district heating network instead. The parameters of the units such as production costs and technical data are given in Table \ref{TechnicalDataandcosts}. The data presented is taken from the energyPRO library, which is based on actualized data from heat-only and CHP units by \cite{technolo56:online}. The electricity prices for 2016 in the DK1 Nordpool area are taken from the Energinet datahub \cite{Dataabou83:online}. The heat demand for 2016 is set to 37500 MWh-heat, comparable to 1 $\text{km}^{2}$ of an urban area with a heat density of 120 TJ/$\text{km}^{2}$.  This amount represents around 2000 family households where 40\% of the heat delivered is domestic hot water and grid losses, which are not weather dependent. The remaining 60\% of the heat delivered is used for space heating which is dependent on the outside temperature.  

We will study two different operation modes of the CHP units. First, \textit{full load} is a restrictive mode in which if the CHP unit is on, it is forced to produced at its maximum capacity ($\underline{P}_{u}=\overline{P}_{u}=2.5$). Many district heating providers in Denmark operate their CHP units only at full load, where it has the highest efficiency. In the other configuration, named \textit{partial load} configuration, both CHP units can produce at partial load between $\underline{P}_{u}$ and $\overline{P}_{u}$. Note that the only change to the method that has to be made to switch between those two modes is setting parameter $\underline{P}_{u}$ equal to $\overline{P}_{u}$.

\begin{table}
	\centering
	\caption{Characteristics of the production units and the thermal storage}
	\begin{tabular}{crrrrrrrrrrrr}
		\toprule
		Unit	&	$C_{u}^{}$	&$\overline{Q}_{u}$	& $\overline{P}_{u}$	&	 $\varphi_{u}$	&	$\overline{S}$	&	 	$\underline{S}$	& $S^{\text{F}}$ & $A^{\text{DH}}_{u}$ &$A^{\text{S}}_{u}$ \\ \hline
		CHP1	&	610.84	&	2.95	&	2.5	&	1.18	&	- &	- &		- & 0  & 1\\
		CHP2	&	610.84	&	2.95	&	2.5	&	1.18	&	- &	- &		- & 0 & 1 \\ 
		GB	    &	404.02	&	19  	&	-	&	-	&	-    &	- & 	- & 1 & 0 \\
		WCB  	&	211.45	&	0.95	&	-	&	-	&	-    &	- &  	- & 0 & 1  \\
		TS  	&	-    	&	-   	&	-	&	-	&	46.93 & 0 & 	46.93	& - & -\\
		\bottomrule
	\end{tabular}  
	\label{TechnicalDataandcosts}
\end{table}

\section{Illustrative example}\label{sec:illustrativeexample}
This example demonstrates our bidding method HURB for the case study described in the previous section. We consider one specific day in March 2016 to illustrate the steps of Algorithm \ref{alg:biddingmethod}. We use the \textit{full load} operation mode and an initial storage level of $s_{0}=10$ MWh-heat.
\begin{figure}
	\centering
	\includegraphics[width=0.65\textwidth]{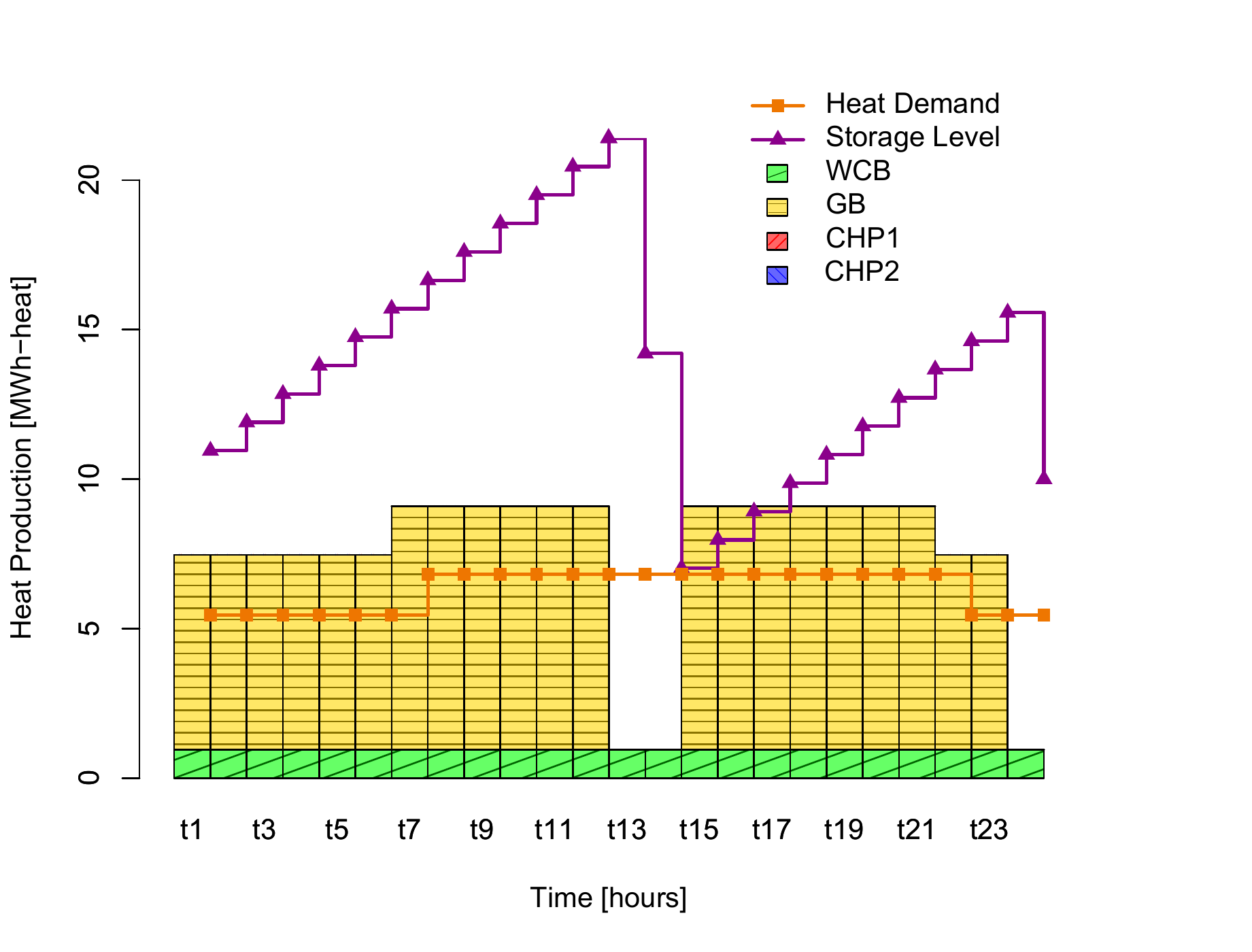}
	\caption{Heat production with no minimum heat restriction and considering no electricity market participation.  }
	\label{FirstSolvingAlgorithm}
\end{figure}

In the initialization of the algorithm, we solve the model without forcing any minimum heat production to the CHP and heat-only units and without considering the electricity price forecast. The optimal solution is depicted in Figure \ref{FirstSolvingAlgorithm} showing that the cheaper WCB and GB are used for production. This obtained solution ($q^{*}_{u,t}$) determines the heat production lower bound ($\underline{Q}_{u,t}$) for future steps in the algorithm. 
The algorithm continues by identifying the heat-only production units and sorting them in descending order of production costs. In our case this results in set $\mathcal{H} = \lbrace \text{GB, WCB} \rbrace$ because the GB has higher production cost (404.62) than the WCB (211.45).

\begin{figure}
	\begin{subfigure}[t]{.51\textwidth}
		\centering
		\includegraphics[width=1\textwidth]{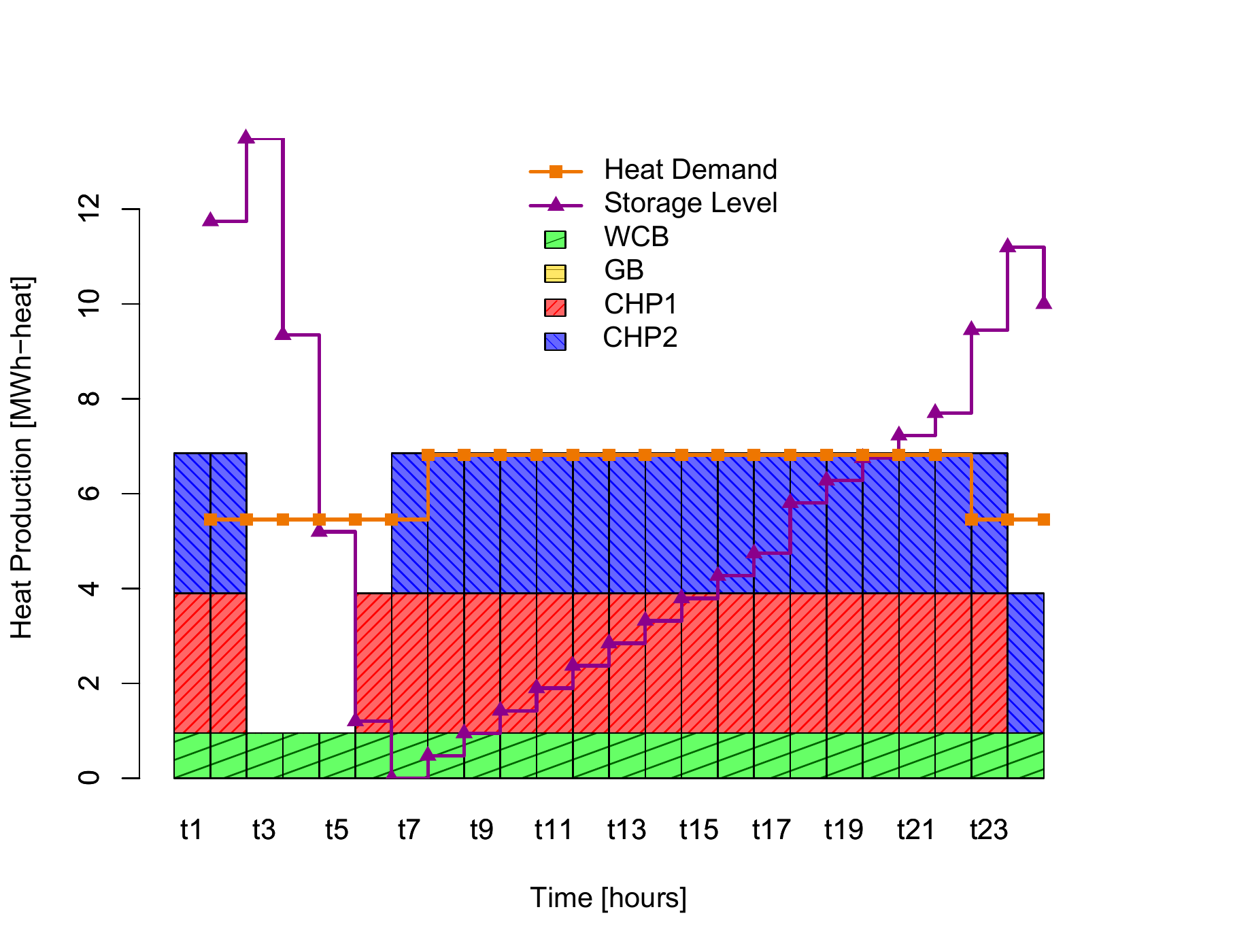}
		\captionsetup{width=.8\textwidth}
		\caption{Heat production and storage level of the system compared to the heat demand}
		\label{fig:firstiterHP}
	\end{subfigure}
	\begin{subfigure}[t]{.51\textwidth}
		\centering
		\includegraphics[width=1\textwidth]{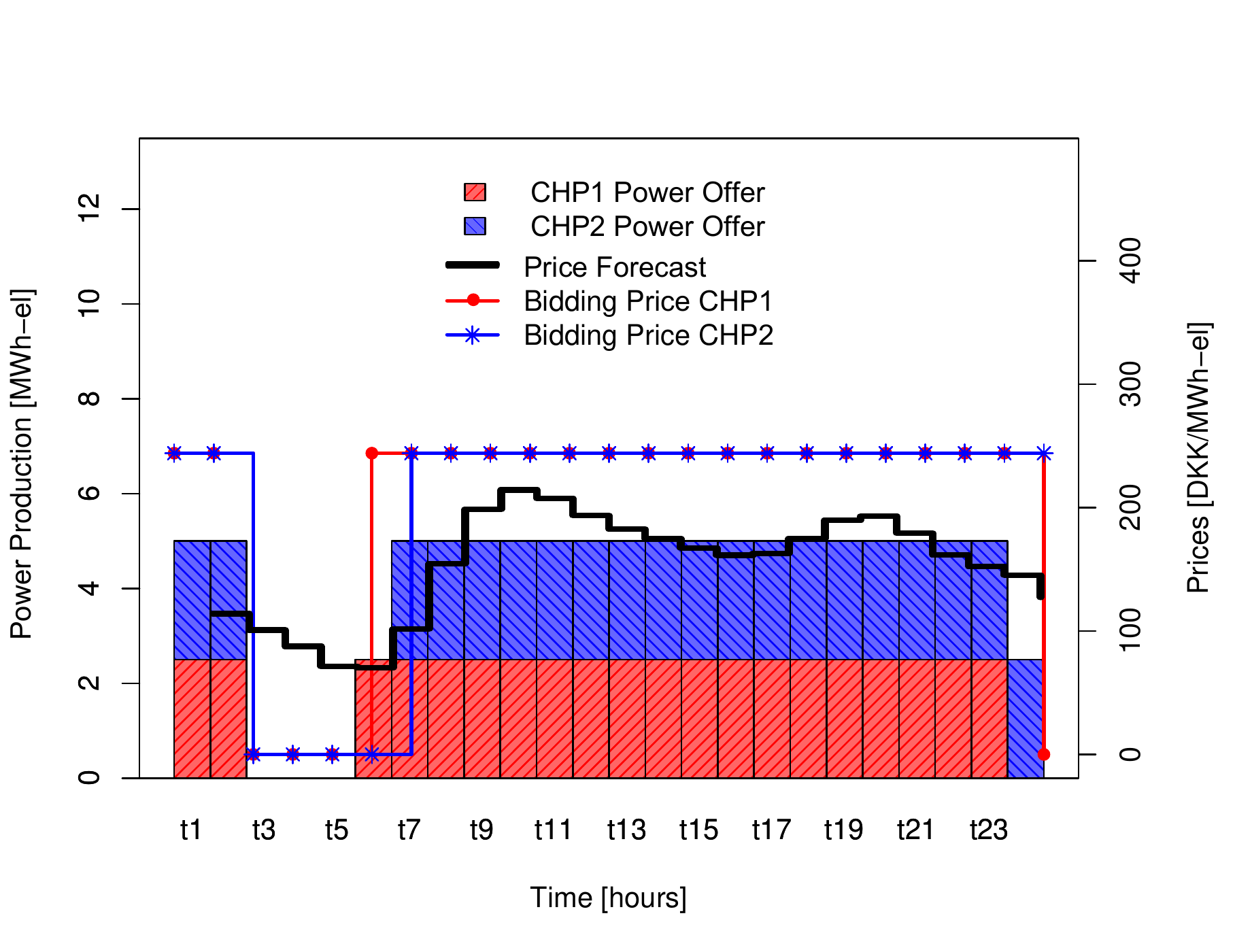}
		\captionsetup{width=.8\textwidth}
		\caption{Power offered and bidding prices for both CHP1 and CHP2 vs. forecast electricity price}
		\label{fig:firstiterPP}
	\end{subfigure}
	\caption{Production data of the system for 24 hours when the gas boiler unit is removed (Iteration 1)} \label{fig:firstiterfigure}
\end{figure}

In the first iteration $i=1$ of the algorithm, we remove the GB from $\mathcal{U}^{\text{H}}$ and solve the optimization problem  \eqref{objectivefunction1}-\eqref{eq2:Heatbalance} for the remaining set of units $\mathcal{U_{}}=\{\text{WCB, CHP1, CHP2}\}$ formed by the union of the sets $\mathcal{U}^{\text{CHP}}=\{\text{CHP1, CHP2}\}$ and $\mathcal{U}^{\text{H}}_1 = \lbrace \text{WCB} \rbrace$, where the minimum production of the WCB is limited to the production from Figure \ref{FirstSolvingAlgorithm}. The resulting optimal heat production and storage operation is shown in Figure \ref{fig:firstiterHP} for the given heat demand.  Based on this information we determine the offers for the day-ahead market. The bidding hours and power amounts for the offers can be directly deferred from the optimal production. The optimal power production, bidding prices and price forecast for this day are depicted in Figure
\ref{fig:firstiterPP}, where both CHP units are used in most of the hours. Except in periods 3 to 5 where both units are off, hour 6 where only CHP1 is on and hour 24 where only CHP2 is on. This means we have to determine offers for CHP1 in hours $\lbrace 1,2,6, \ldots, 23\rbrace$ and for CHP2 in hours $\lbrace 1,2,7, \ldots, 24\rbrace$. The power amount is always 2.5 MWh per CHP. The bidding price for the two CHP units is the same as they have the same operational cost and determined by subtracting the cost of the removed GB from the CHP operational cost and adjusted from price per MWh heat to price per MWh power. This results in the following bidding price \eqref{biddingproce1stiter} (see also Figure \ref{fig:firstiterPP}).
\begin{align}
\forall u \in \mathcal{U}^\text{CHP}\quad 
b_{u} &=\big(C_{u}-C^{}_{\text{GB}}\big)\cdot \varphi_{u}\label{biddingproce1stiter}\\ 
&= (610.84 - 404.02) \cdot 1.18 = 244.045 \nonumber
\end{align}
The iteration finishes with a total number of 40 hourly offers where each offer has a power amount of 2.5 MWh and a price of 244.045 DKK.

\begin{figure}
	\begin{subfigure}[t]{.51\textwidth}
		\centering
		\includegraphics[width=1\textwidth]{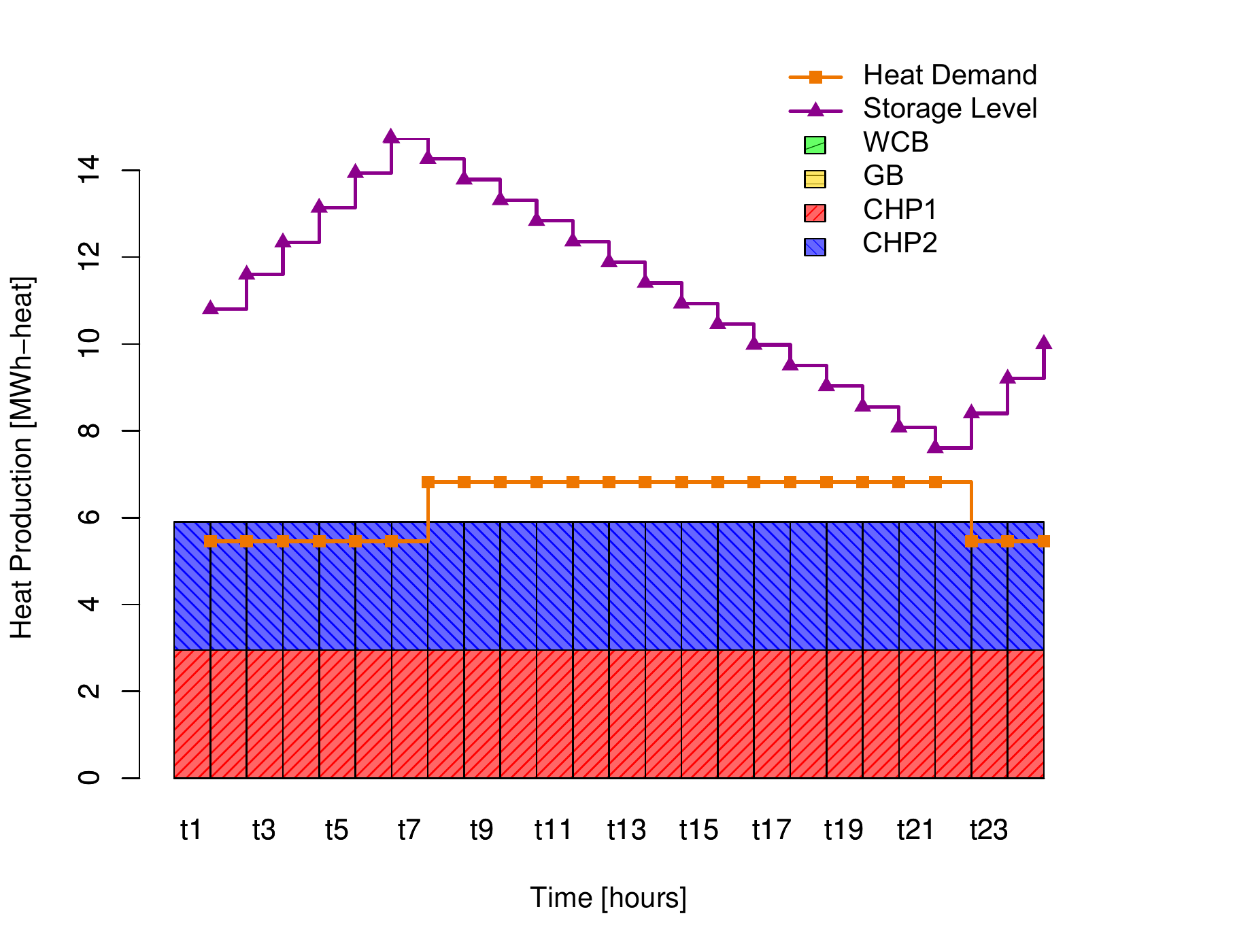}
		\captionsetup{width=.8\textwidth}
		\caption{Heat production and storage level of the system compared to the heat demand}
		\label{fig:seconditerHP}
	\end{subfigure}
	\begin{subfigure}[t]{.51\textwidth}
		\centering
		\includegraphics[width=1\textwidth]{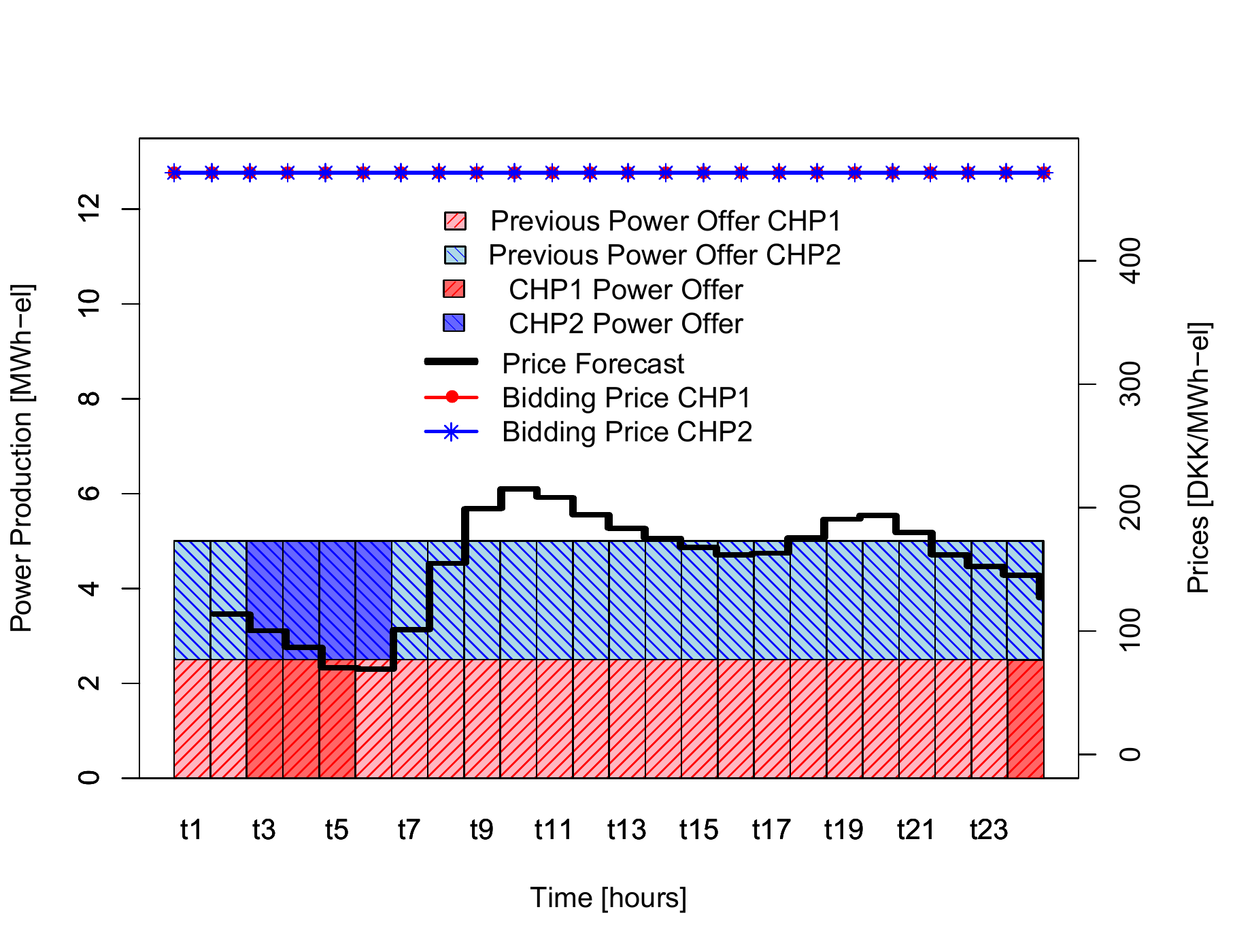}
		\captionsetup{width=.8\textwidth}
		\caption{Power offered and bidding prices for both CHP1 and CHP2 vs. forecast electricity price}
		\label{fig:seconditerPP}
	\end{subfigure}
	\caption{Production data of the system for 24 hours when both the gas boiler unit and the wood chip boiler are removed (Iteration 2)} \label{fig:seconditerfigure}
\end{figure}

In the second iteration $i=2$, we remove the next most expensive heat-only unit, which is the WCB, reducing the set $\mathcal{U}^{\text{H}}_2$ to an empty set. Thus, in this iteration we solve the optimization model \eqref{objectivefunction1}-\eqref{eq2:Heatbalance} using just the CHP units $\mathcal{U} = \mathcal{U}^{\text{CHP}}$. The resulting heat and power production is shown in Figure \ref{fig:seconditerfigure}. The hours and power amounts for the new offers are the difference between the power production in the previous iteration and this iteration. This means we add new offers for CHP1 in hours $\lbrace 3,4,5,24 \rbrace$ and for CHP2 in hours $\lbrace 3,4,5,6 \rbrace$. All new offers will have the amount 2.5 MWh, because the respective CHP was off before and is now producing at full load. The bidding price is calculated in equation \eqref{biddingproce2stiter} and this time based on the cost of the WCB.  
\begin{align}
\forall u \in \mathcal{U}^\text{CHP}\quad 
b_{u} &=\big(C_{u}-C^{}_{\text{WCB}}\big)\cdot \varphi_{u}\label{biddingproce2stiter}\\ 
&= (610.84 - 211.45) \cdot 1.18 = 471.279 \nonumber
\end{align}
Thus, the second and last iteration adds eight new offers to the set.

\begin{figure}
	\begin{subfigure}[t]{.54\textwidth}
		\centering
		\includegraphics[width=1\textwidth]{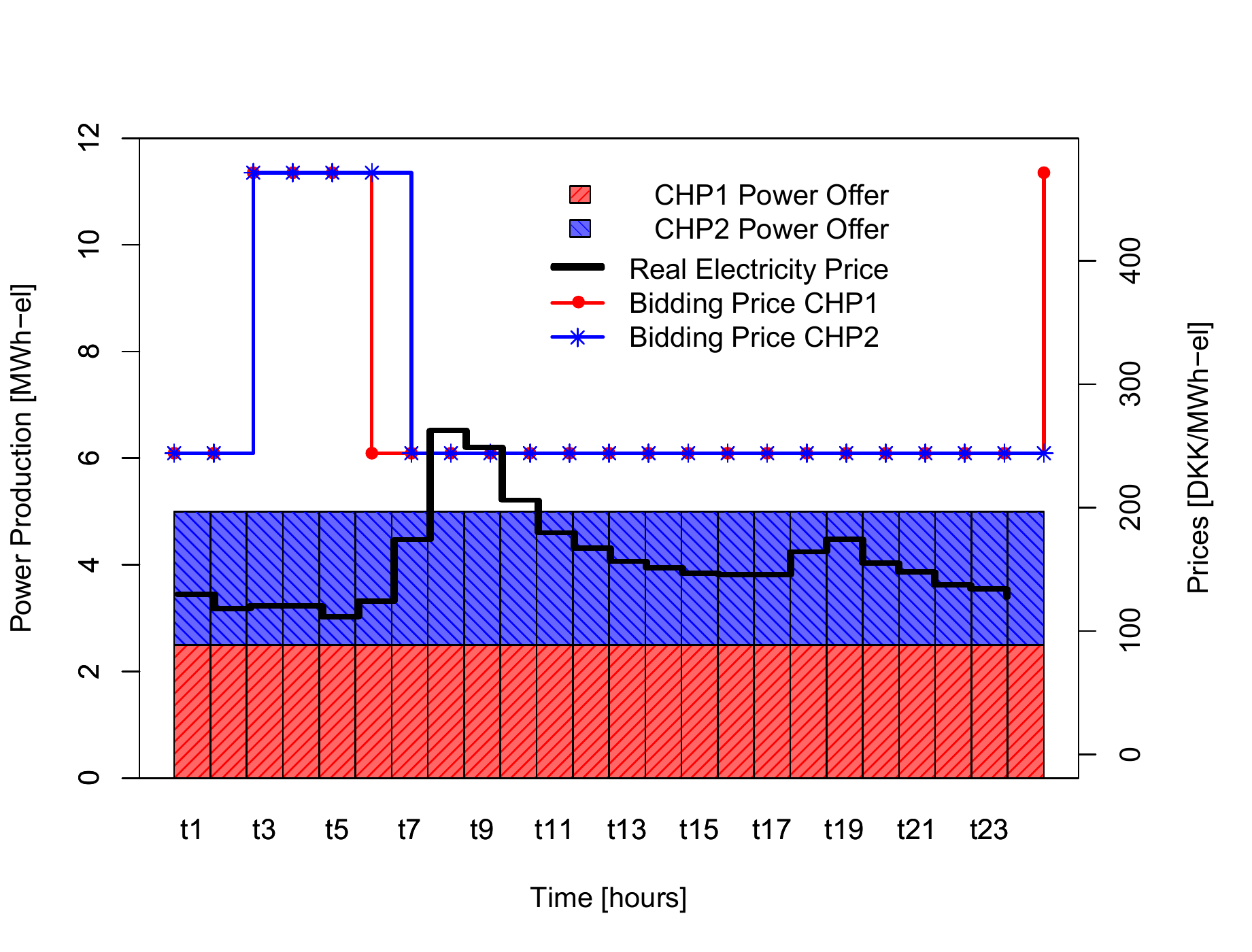}
		\captionsetup{width=.7\textwidth}
		\caption{Power offered and bidding prices for both CHP1 and CHP2 vs. real electricity prices}
		\label{fig:totalPP}
	\end{subfigure}
	\begin{subfigure}[t]{.54\textwidth}
		\centering
		\includegraphics[width=1\textwidth]{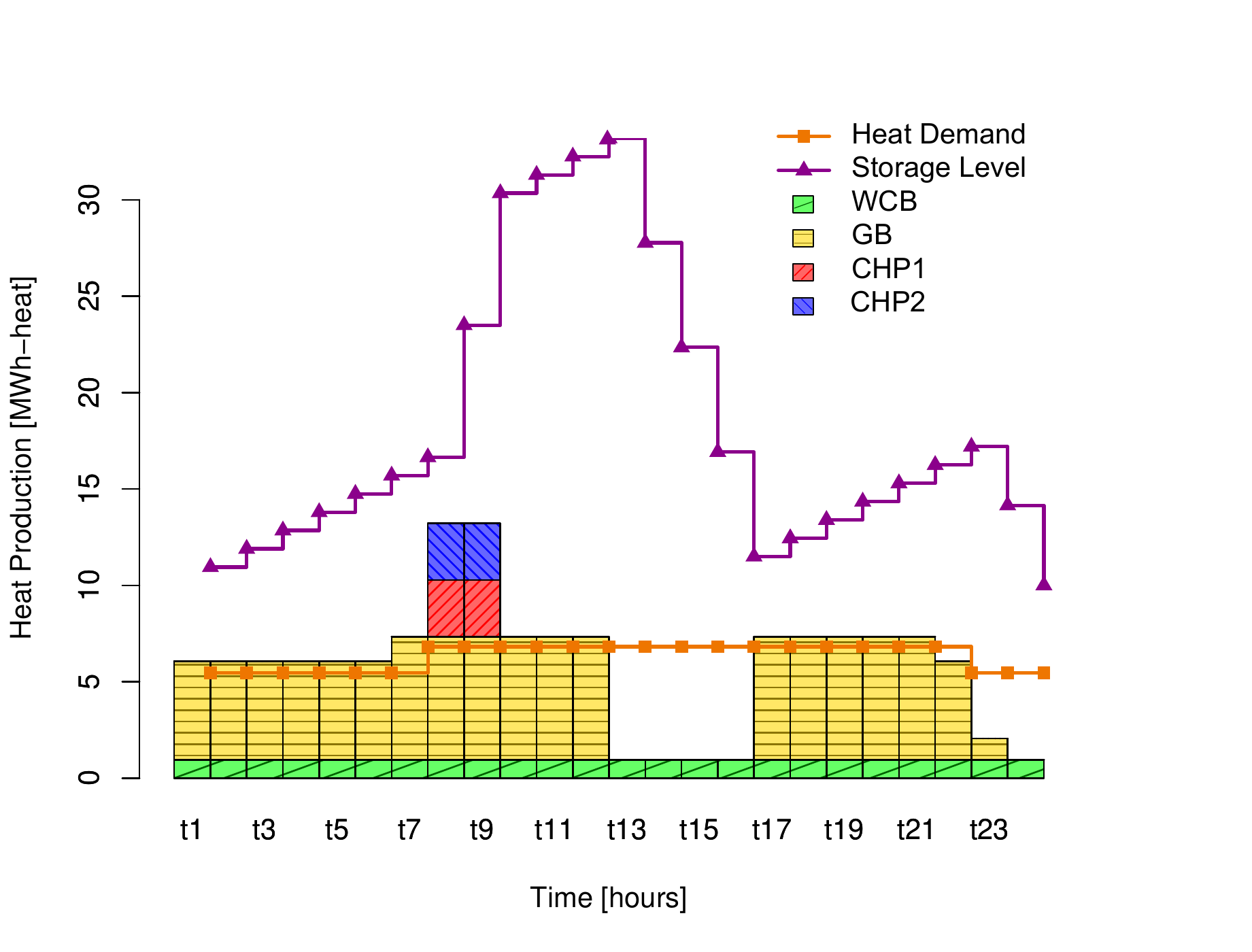}
		\captionsetup{width=.7\textwidth}
		\caption{Real heat production and storage level of the system compared to the heat demand}
		\label{fig:totalHP}
	\end{subfigure}
	
	\caption{Final production of the system using real data for the day studied.} \label{fig:totalfigure}
\end{figure}

In total, the bidding method resulted in 48 bids. Figure \ref{fig:totalPP} shows the final offers and bidding prices presented to the market. The price depicted is now the real realization instead of the forecast used in the optimization. For the operation of the system this implies that if the real price is greater or equal than the bidding price for an specific period of time, the unit is committed to produce (in this case for hours 8 and 9). In the other periods, we are free to determine the operation of the system. Based on the optimization of model \eqref{objectivefunction1}-\eqref{eq2:Heatbalance} with fixed power production in periods 8 and 9, we obtain the heat production schedule in Figure \ref{fig:totalHP}. Except in the committed hours of CHP production, we use the cheaper WCB and GB in the other periods to cover the heat demand. Furthermore, we can see from the storage level also shown in Figure \ref{fig:totalHP} that the thermal storage is used to avoid more expensive GB production in some periods.

\section{Numerical Results}\label{sec:results}

In this section, we analyze the performance of our bidding method and compare it to the following proposed methods from literature:

\begin{itemize}
	\item \textit{Method A:}  Bidding strategy based on the uncertainty bound given by forecasted prices \cite{conejo2002price}.
	\item \textit{Method B:}  Use of confidence intervals on price forecast to create bidding curves \cite{rodriguez2004bidding}.
	\item \textit{Method C:} Forecasted electricity price as bidding prices \cite{schulz2016optimal}. 
	\item \textit{Method D:}  Use of price scenario intervals to create bidding curves \cite{dimoulkas2014constructing}.
	\item \textit{Method E:} Piece-wise linear function of costs used as bidding prices \cite{ravn2004modelling}.  
\end{itemize}

Although methods A to D are developed and designed for standalone CHP or thermal production units, they could be used also by district heating operators to determine the bids for the CHP units without taking the other units into account. We compare the HURB method to these methods to show that district heating operators have a benefit from using a method that takes the  entire portfolio of units into account and should therefore use more specialized methods instead of the general applicable bidding methods. However, the HURB method would not be applicable to the general case without other heat producing units as it is considered in \cite{conejo2002price,rodriguez2004bidding,schulz2016optimal,dimoulkas2014constructing}.

All methods base their offers on the optimal operation of the system. For a fair comparison of the bidding strategy, we use the model \eqref{objectivefunction1}-\eqref{eq2:Heatbalance} to determine the optimal power production for all methods. In this way, we focus on just comparing the creation of bids and not the underlying operational planning problem. For the above mentioned methods A to E, no iterative procedure is proposed to replace heat-only units. Furthermore, the bidding prices are determined in different ways. 
For further information about the details of the methods we refer the reader to the mentioned references.

The following comparison of the methods is based on the totals costs by using the bidding methods for each day of the planning horizon. This daily optimization can be seen individually by just considering 24 time periods. However, we use a receding approach and, thus, include more time periods in the operational optimization to consider the operation of the thermal storage over several days. In the end, only the offers for the first 24 hours are send to the market. For a fair comparison we ensure that the storage level at the end of the horizon is the same for all methods. 

The daily usage of the bidding method also allows us to update the forecasts for electricity prices every day. For this experiment, we use time series analysis to create a SARIMA model with one day seasonality following the same approach as in \cite{nogales2002forecasting} including harmonic external regressors in the form of Fourier series \cite{weron2007modeling} to describe weekly seasonality ($T$=168 hours). The model to predict electricity prices $\lambda_t$ is 
\begin{align}
\lambda_{t}  & =  \mu+\phi_{1}\lambda_{t-1}+\phi_{2}\lambda_{t-2}+\phi_{24}\lambda_{t-24}+\theta_{1}\varepsilon_{t-1}+\theta_{2}\varepsilon_{t-2}+\theta_{24}\varepsilon_{t-24} \notag\\
& +\sum_{k=1}^{K}\bigg[\alpha_{k}\sin{\bigg(\frac{2\pi k t}{168}\bigg)}+\beta_{k}\cos{\bigg(\frac{2\pi k t}{168}}\bigg)\bigg]
\label{eq:Fourierseries}
\end{align}
where $\lambda_{t}$ is the estimated electricity price for period $t$. The coefficient $\mu$ is the intercept and $\lambda_{t-1}$, $\lambda_{t-2}$, $\lambda_{t-24}$, $\varepsilon_{t-1}$, $\varepsilon_{t-2}$ and $\varepsilon_{t-24}$ correspond to the autoregressive and moving average terms, respectively. $K$ is a natural number chosen by minimizing the Akaike information criterion (AIC) and determines the number of Fourier terms considered. Finally, $\alpha_{k},\beta_{k},\phi_{1},\phi_{2},\phi_{24},\theta_{1},\theta_{2}$ and $\theta_{24}$ are the forecast parameters.  More sophisticated forecasting methods based on temporal dependencies such as probabilistic forecasting can be used. We refer the reader to, e.g., \cite{jonsson2014predictive} where the authors create a model to predict day-ahead electricity prices based on reliable probability density forecasts.

Based on this, the evaluation process can be described as follows:
\begin{enumerate}
	\item Choose length of receding horizon $|\mathcal{T}|$ (a multiplier of 24 to cover entire days).
	\item Update forecast parameters of electricity prices for day $d$ using the most recent observations and predict for $|\mathcal{T}|$ periods. \label{update}
	\item Apply bidding method (either HURB or method A-E) for $|\mathcal{T}|$ periods. \label{solving}
	\item Evaluate the set of offers $\mathcal{O}$ for real prices of day $d$.
	\item Get the system cost for day $d$.
	\item Update $d=d+1$ and go to step \ref{update}. 
\end{enumerate}

\subsection{Evaluation for the year 2016}

We first apply the methods HURB and A to E for the year 2016. All methods are tested with different lengths of receding horizon, namely 1 to 15 days, and use the mentioned SARIMA model for forecasting the electricity prices. 

\begin{figure}
	\begin{subfigure}{.55\textwidth}
		\centering
		\includegraphics[width=1.0\textwidth]{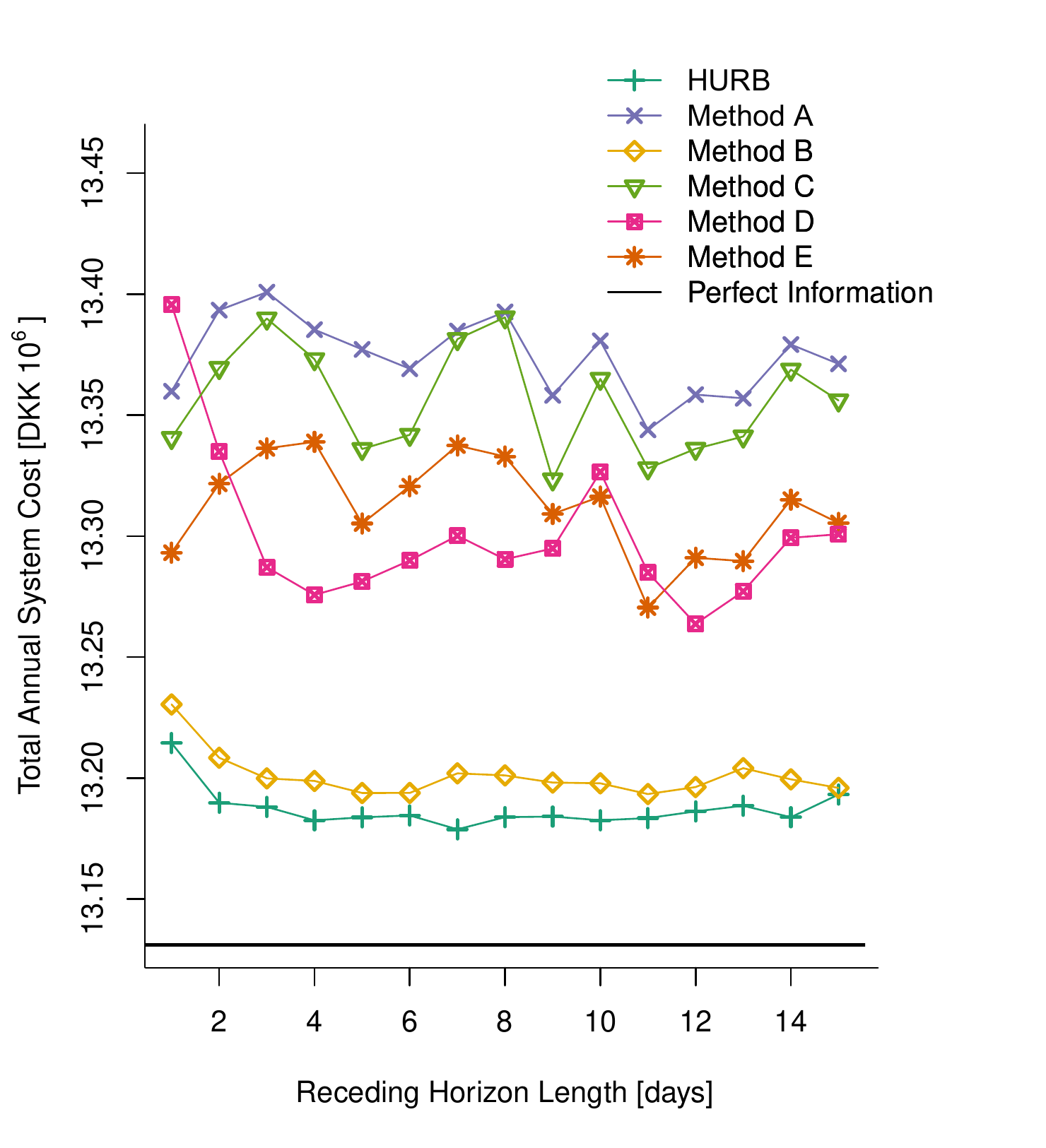}
		\captionsetup{width=.8\textwidth}
		\caption{\textit{Full load configuration} }
		\label{fig:year_full}
	\end{subfigure}
	\begin{subfigure}{.55\textwidth}
		\centering
		\includegraphics[width=1.0\textwidth]{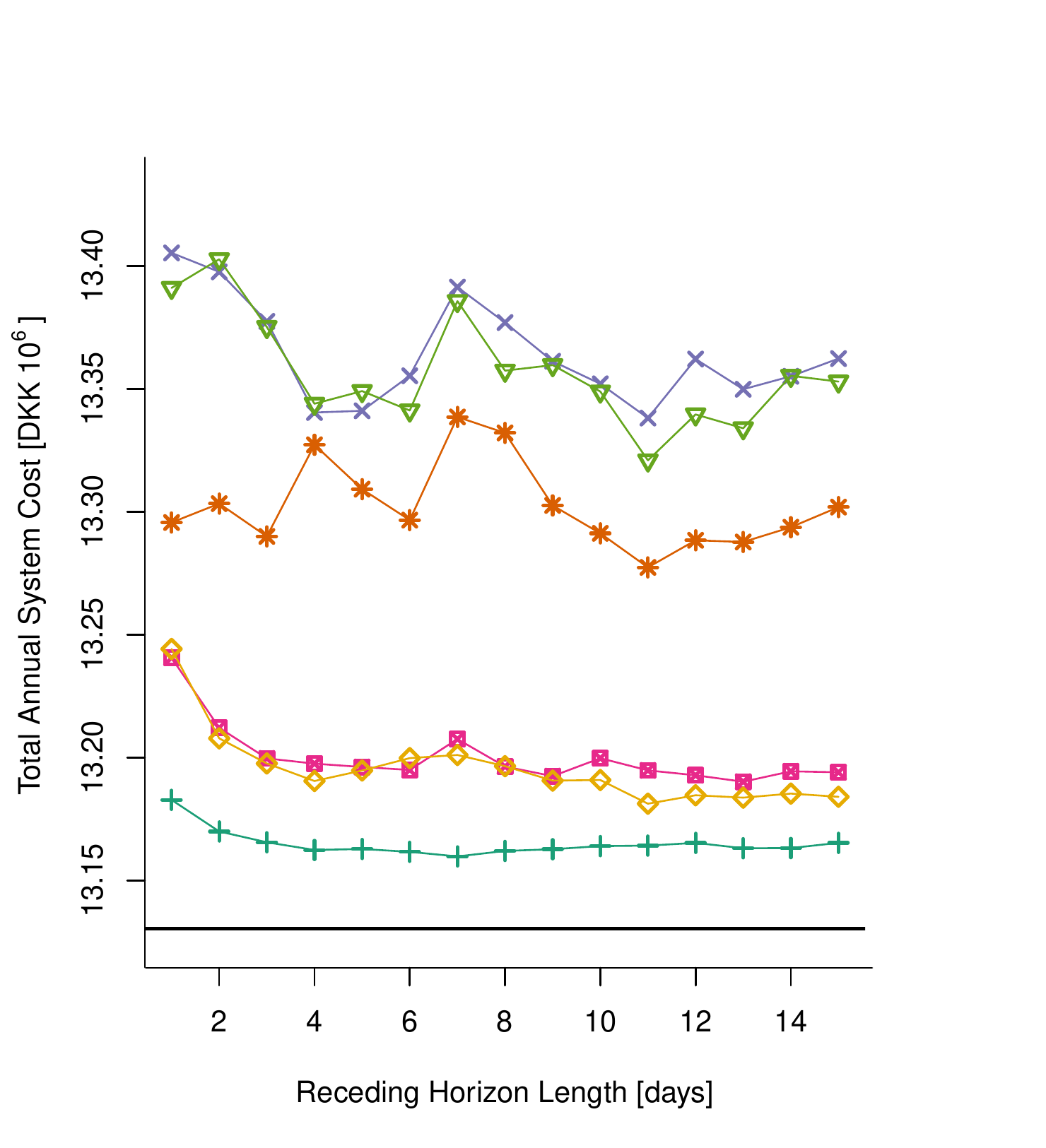}
		\captionsetup{width=.8\textwidth}
		\caption{\textit{Partial load configuration}}
		\label{fig:year_partial}
	\end{subfigure}
	\caption{Annual system cost per method and length of receding horizon (1 to 15 days). The black line shows the cost with perfect information (solving the entire year at once with real electricity prices). } \label{fig:2016}
\end{figure}

Figure \ref{fig:2016} shows the annual system costs for all methods. Furthermore, the figure distinguishes between the two configurations \textit{full load} (\ref{fig:frontier1}) and \textit{partial load} (\ref{fig:frontier2}). For all lengths of receding horizon and both configurations, our method (HURB) yields the best result regarding total annual system costs, closely followed by method B. Those two methods are also very close to the total annual costs when having perfect information (black line), which are the theoretical minimal costs of the system when we solve the optimal operation for the entire year at once and have perfect information about electricity prices. Figure \ref{fig:2016} also shows that all methods benefit from using a receding horizon of more than one day, i.e., the solution can be improved by already considering future days for the operation of the storage.

\subsection{Evaluation of further electricity price sets}

To get a more general result of the comparison, we evaluate a more diverse set of electricity prices, i.e., we use different samples of realization of electricity prices to evaluate how the method would perform if the electricity prices would have turned out differently. Therefore, we repeat our experiment again for the month January 2016 using the same price forecasting method as before. We use 144 samples for realizations of electricity prices obtained from different markets in the Nordpool region to evaluate the performance of the methods in the market. We use data from the regions Sweden (SE4) and Germany (DE) for the months in 2013 to 2015 and 2017. As the cost per month differ based on the prices, we compare the methods using the difference between the cost obtained by the method and the cost obtained by using perfect information regarding electricity prices (solving the entire month at once with the respective electricity prices instead of using a forecast). 

\begin{figure}
	\begin{subfigure}{.515\textwidth}
		\centering
		\includegraphics[width=1.0\textwidth]{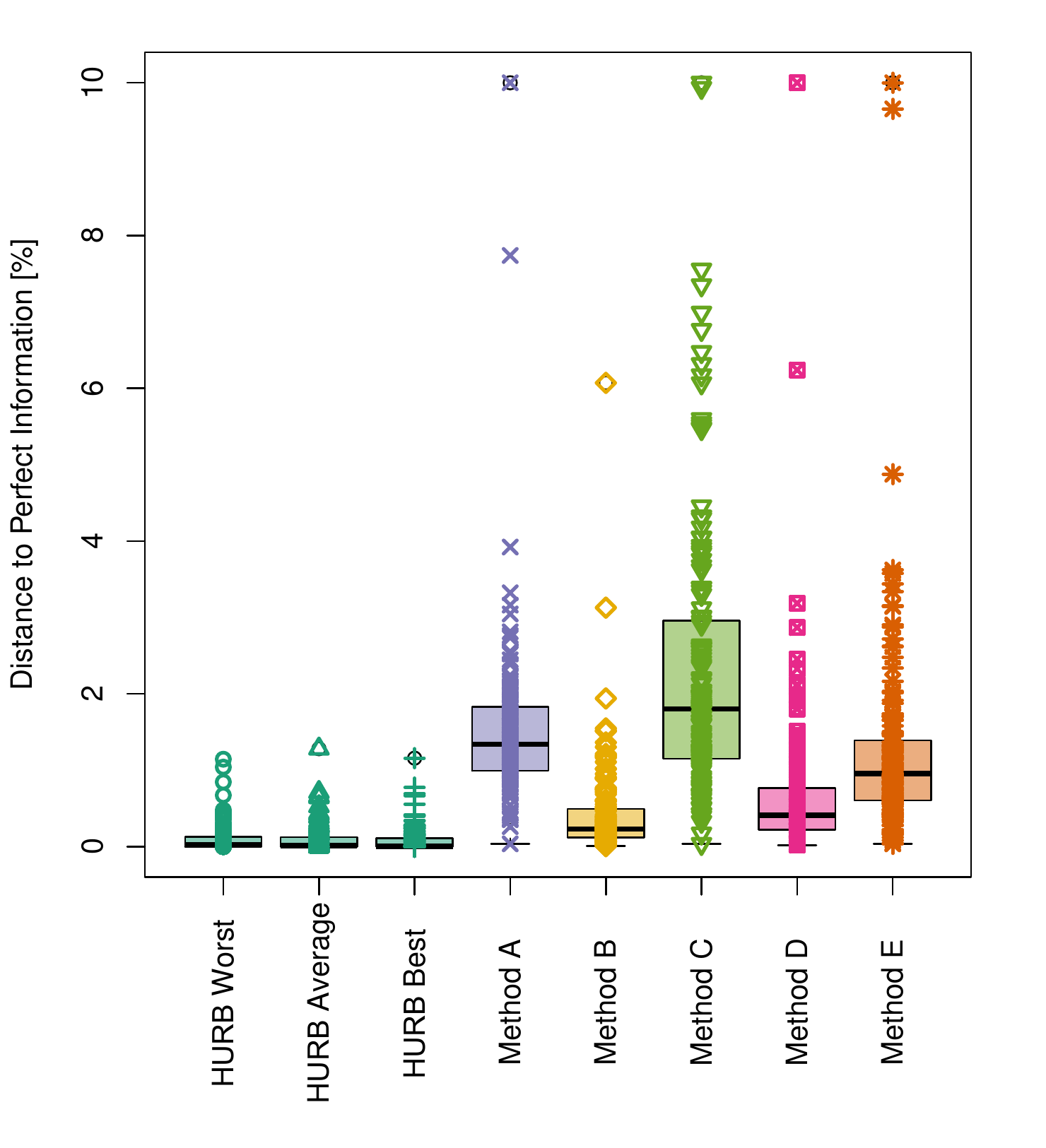}
		\captionsetup{width=.8\textwidth}
		\caption{\textit{Full load configuration} }
		\label{fig:frontier1}
	\end{subfigure}
	\begin{subfigure}{.515\textwidth}
		\centering
		\includegraphics[width=1.0\textwidth]{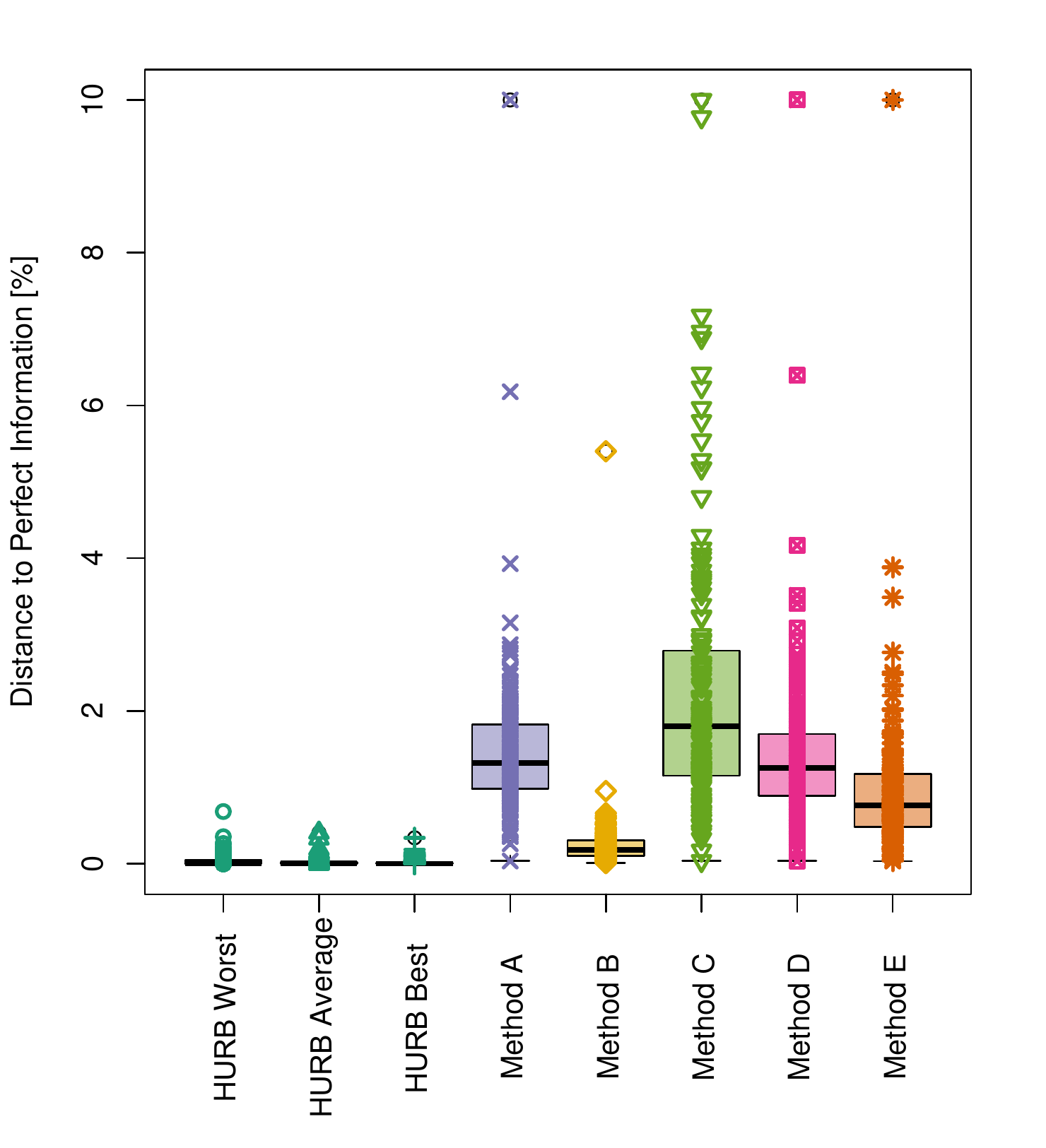}
		\captionsetup{width=.8\textwidth}
		\caption{\textit{Partial load configuration}}
		\label{fig:frontier2}
	\end{subfigure}
	\caption{Distance in \% of montly system cost per method to the cost obtained by using perfect information. All methods are run with 1 to 15 days of receding horizon. For method A to E and HURB Best, the results show the cost for the best receding horizon. For HURB Worst, the worst receding horizon is chosen and HURB Average shows the mean difference to the perfect information solution over all lengths of receding horizon.} \label{fig:frontier}
\end{figure}

Figure \ref{fig:frontier} shows the distance in \% to the cost using perfect information. All methods were again evaluated using different lengths of receding horizon. For method A to E and HURB Best, the results show the cost for the best receding horizon length. For HURB Worst, the worst receding horizon is chosen and HURB Average shows the mean difference to the perfect information solution over all lengths of receding horizon. The specific number of days are also mentioned in Table \ref{tab:results_mode}. Furthermore, the figure distinguishes between the two configurations \textit{full load} (\ref{fig:frontier1}) and \textit{partial load} (\ref{fig:frontier2}).

For both configurations HURB yields the best result, i.e., the cost are closest to perfect information cost. Even with the worst configuration of receding horizon length (HURB Worst), it leads to better results for most of the samples compared to the other methods. Our method is closely followed by method B, which outperforms the other methods A and C to E in most of the samples. 

\subsection{Discussion}
\begin{table}
	\footnotesize
	\begin{adjustbox}{max width=\textwidth}
		\begin{tabular}{lrrrrrrrrrr}
			\toprule
			Method &                    \multicolumn{ 5}{c}{Full load} &                 \multicolumn{ 5}{c}{Partial load} \\
			\cmidrule(r){2-6}\cmidrule(r){7-11}
			
			& RH & \multicolumn{ 2}{c}{CHP 1} & \multicolumn{ 2}{c}{CHP 2} & RH & \multicolumn{ 2}{c}{CHP 1} & \multicolumn{ 2}{c}{CHP 2} \\
			\cmidrule(r){3-4}\cmidrule(r){5-6}\cmidrule(r){8-9}\cmidrule(r){10-11}
			& &     Offers &        Won &     Offers &        Won &   &  Offers &        Won &     Offers &        Won \\\midrule
			
			HURB Worst &  1 &  98.44 & 41.88 & 98.41 & 41.81 & 1 & 98.91 & 41.95&98.70 & 41.91 \\
			
			HURB Avg. &  - &   99.64  &  42.16  &  99.62  &  42.07 & - &  99.79  &  42.19  &  99.75  &  42.15 \\
			
			HURB Best &  10 &    99.83  &  42.27  &  99.82  &  42.19  & 10 &  99.89  &  42.28  &  99.87  &  42.26 \\
			
			Method A & 11 &      44.88  &  39.29  &  44.87  &  39.27  & 10 &  44.92  &  39.34  &  44.92  &  39.31 \\
			
			Method B & 14 &      82.25  &  35.25  &  82.15  &  35.25  & 5 & 82.52  &  35.85  &  82.40  &  35.82 \\
			
			Method C &  2 &      44.60  &  18.33  &  44.66  &  18.36  & 12 &  45.02  &  18.54  &  45.01  &  18.53 \\
			
			Method D & 10 &      74.62  &  35.61  &  74.61  &  35.60  & 12 & 75.55  &  26.56  &  75.55  &  26.55 \\
			
			Method E & 5 &      44.62  &  31.87  &  44.63  &  31.89  & 5 & 44.84  &  32.58  &  44.83  &  32.57 \\
			\bottomrule
		\end{tabular}  
	\end{adjustbox}
	
	\caption{Percentage of hours with offers and won bids in the investigated month averaged over all samples (RH = receding horizon).}
	\label{tab:results_mode}
\end{table}

To explain why our methods outperforms the others, we state the average number of hours and won bids for each method over all samples in Table \ref{tab:results_mode}. The values are given in percentage of hours per month. 

These values show that HURB  submits more offers to the market than the other methods. On average HURB places offers in more than 98\% of the hours, while the other methods have fewer bids. In particular, method A, C and E only place offers in less than 50\% of the hours. HURB also wins more bids in the market compared to the others (in more than 40\% of the hours). Method B \cite{rodriguez2004bidding}, which is the second best regarding cost, wins on average in less than 36\%  of the hours. Thus, HURB achieves a higher income from the market. Although, method A \cite{conejo2002price} also yields won bids close to 40\% of the hours like HURB, these bids are not always profitable for the district heating operator. As method A was developed for a standalone thermal unit it does not take the other units into account to calculate the overall system costs. The bidding prices are set based on the confidence interval of the electricity price forecast, which could lie below the marginal cost of the system. Therefore, method A is not advisable for district heating operators.  The same  holds for method C \cite{schulz2016optimal}, which uses the forecast as bidding price. 

Method B  and D \cite{dimoulkas2014constructing} are using the forecast on electricity prices to determine when and at which price it is profitable to produce with the CHP units. While method B uses intervals of the probability density function of the price forecast, method D uses sampled scenarios. They analyze the results for the intervals/scenarios and present an offer with the respective bidding price of the interval/scenario, if the CHP units produces according to their optimization method. As these methods also do not make use of the other units in the district heating system, they offer only in hours where the forecast indicates profits. In the HURB method, we explicitly consider the district heating system and make use of the fact that we have to produce the heat either way. Therefore, the electricity price forecast does not determine the amount or price of the offers but only the hours in which to place them.  This is also the explanation, why the HURB method bids in more hours than method E \cite{ravn2004modelling}. Although method E also uses the marginal cost as bidding prices and makes uses of the district heating system, the creation of offers regarding bidding amounts and offers is different. Method E places offers at the marginal prices of the units, whenever it is profitable with respect to the price forecast. In contrast, HURB offers as much electricity which is needed to replace the entire heat demand by CHP production, which results in more offers in total. Furthermore, the prices set by HURB ensure no losses regarding the cost minimal production without trading, but only profits if bids are won.

\section{Summary and outlook}
\label{conclusion}
In this work, we propose a new day-ahead market bidding strategy, named Heat Units Replacement Bidding (HURB) method, for CHP units that are operated jointly with other heat or heat-only production units. Our bidding method is based on replacing heat production from heat-only units through CHP production in a iterative manner. The method uses a mathematical program to determine the optimal operation of the portfolio of units. In order to evaluate our bidding method, we implement other bidding methods proposed in  literature and compare them in a common case study. We state the results for an entire year by executing the bidding method on a daily basis with different lengths of receding horizon. Furthermore, we perform an out-of-sample test to get more general results. The results show that compared to the other bidding methods, our method yields the lowest cost and most won bids. 

To extend the use of this bidding method, we propose three future research directions. First, the consideration of the balancing market offers further opportunities to reduce the operational cost. Second, the modelling of block bids, start-up cost and minimum operation times for the CHP units is a valuable extension to cover more instruments on the electricity  markets. Finally, the increasing presence of solar thermal units at district heating operators introduces an additional source of uncertainty regarding the production amounts. This production has to be taken into account appropriately while determining the offers to the market. Methods based on stochastic programming are able to incorporate this uncertainty by including potential scenarios.

\section*{Acknowledgments}
 This work is funded by Innovation Fund Denmark through the CITIES research center (no. 1035-00027B).


\bibliography{bibliography}

\end{document}